\input amstex
\documentstyle{amsppt}
\input bull-ppt
\keyedby{bull334e/pah}
\topmatter
\cvol{28}
\cvolyear{1993}
\cmonth{January}
\cyear{1993}
\cvolno{1}
\shorttitle{Factorizations of invertible operators}
\cpgs{75-83}
\ratitle
\title Factorizations of invertible operators\\
 and $K$-theory of $C^*$-algebras
\endtitle
\author Shuang  Zhang
\endauthor
\address Department of Mathematical Sciences,  University 
of Cincinnati,
Cincinnati, Ohio 45221-0025\endaddress
\ml szhang\@ucbeh.san.uc.edu\endml
\date February 12, 1992. The main results of this article
were presented at the AMS meeting at Springfield, 
Missouri, March
27--28, 1992\enddate
\thanks Partially supported by  NSF\endthanks
\abstract Let $\Scr A$ be a unital C*-algebra. 
We describe  \it K-skeleton factorizations  \rm 
of all  invertible operators on a  Hilbert C*-module $\Scr 
H_{\Scr A}$, 
in particular on $\Scr H=l^2$, 
with the Fredholm index as an invariant. We then outline 
the isomorphisms
$K_0(\Scr A) \cong \pi _{2k}([p]_0)\cong \pi _{2k}
({GL}^p_r(\Scr A))$ 
and $K_1(\Scr A)\cong \pi _{2k+1}([p]_0)\cong 
\pi _{2k+1}(GL^p_r(\Scr A))$ for $k\ge 0 $,
where $[p]_0$ denotes the class of all compact 
perturbations of 
a projection $p$ in the infinite Grassmann space 
${Gr}^{\infty }(\Scr A)$ 
and $GL^p_r(\Scr A)$ stands for the group of all those 
invertible 
operators on $\Scr H_{\Scr A}$ essentially commuting with 
$p$.\endabstract
\endtopmatter

\document
\heading 1. Introduction\endheading

Throughout, we assume that $\Scr A$ is any  unital 
C*-algebra.  
Let $\bold {\Scr H_{\Scr A}}$ be the Hilbert (right) $\Scr 
A$-module  
consisting  of all $l^2$-sequences in $\Scr A$; i.e., 
$\Scr H_{\Scr A}:=\{\{a_i\}:\sum_{i=1}^\infty a_i^*a_i \in 
\Scr A\},$ 
on which  an $\Scr A$-valued inner product and a norm are 
naturally 
defined by  $<\{a_i\}, \{b_i\}>=\sum_{i=1}^\infty 
a_i^*b_i$ and
$\Vert \{a_i\}\Vert =\Vert (\sum_{i=1}^\infty 
a_i^*a_i)^{1/2}\Vert$.
Let $\bold {\Scr L(\Scr H_{\Scr A})}$ stand for the 
C*-algebra 
consisting of all bounded  operators on 
$\Scr H_{\Scr A}$  whose  adjoints exist, and let 
$\bold {\Scr K(\Scr H_{\Scr A})}$ denote the closed linear 
span 
of all finite rank operators on $\Scr H_{\Scr A}$, 
respectively.  
In case $\Scr A$ is the algebra $\bold{C}$ of all complex 
numbers, 
$\Scr H_{\Scr A}$ is the separable, infinite-dimensional 
Hilbert space 
$\Scr H=l^2$; correspondingly, $\Scr L(\Scr H_{\Scr A}) $ 
reduces 
to the algebra  $\bold {\Scr L(\Scr H)}$ of all bounded  
operators 
on $\Scr  H$,  and $\Scr K(\Scr H_{\Scr A})$ reduces to the 
algebra $\bold {\Scr K}$ of all compact operators on $\Scr 
H$.   
Each element in  $\Scr L(\Scr H_{\Scr A})$ can be 
identified with an 
infinite, \it bounded \rm matrix whose entries are 
elements in 
$\Scr A$ [Zh4, \S 1]. This identification can be realized 
by C*-algebraic
 techniques and the two important *-isomorphisms
$\Scr L(\Scr H_{\Scr A})\cong M(\Scr A\otimes \Scr K)$ 
and $\Scr K(\Scr H_{\Scr A})\cong \Scr A\otimes \Scr K
\ = (\underset {\rightarrow }\to {\text{lim}}M_n(\Scr 
A))^-;$
where ${M(\Scr A\otimes \Scr K)}$ is the multiplier 
algebra of  
$\Scr A\otimes \Scr K$ [Kas]. 
For more information about multiplier algebras the reader 
is referred
to  [APT, Bl, Cu1, El, Br2, Pe1, OP, L, Zh4--5], among 
others. 
The set of projections 
$${Gr}^{\infty }(\Scr A):=\{p\in \Scr L(\Scr H_{\Scr A}):\
p=p^2=p^*\ \ \text{and}\ \ p\sim 1\sim 1-p\}$$  is called 
the 
\it infinite Grassmann space associated with $\Scr A$; \rm
where $\lq q\sim p$' is the well-known Murray-von Neumann 
equivalence
of two projections; i.e.,   there exists a partial isometry 
 $v\in \Scr L(\Scr H_{\Scr A})$ such that $vv^*=p$ and 
$v^*v=q$.
If $\Scr A=\bold{C}$, then $Gr^{\infty }(\Scr A)$ reduces 
to  the 
well-known Grassmann  space  ${Gr}^{\infty }(\Scr H)$
consisting of all projections on $\Scr H$ with an infinite 
dimension and 
an infinite codimension.

\heading 2. Factorizations and $K$-theory\endheading

Let $p\in Gr^{\infty }(\Scr A)$. If $x$ is any element in 
$\Scr L(\Scr H_{\Scr A})$, with respect to the 
decomposition 
$p\oplus (1-p)=1$    one can write $x$ as a  $2\times 2$ 
matrix, say 
 $ (\smallmatrix a&b\\c&d\endsmallmatrix  ),$ where 
$a=pxp$, $b=px(1-p)$, $c=(1-p)xp$, and $d=(1-p)x(1-p)$.
A unitary operator $u = (\smallmatrix 
a&b\\c&d\endsmallmatrix )$  
is called a \it K-skeleton unitary  along $p$, \rm 
if both $b$ and $c$  are some partial isometries 
in $\Scr A\otimes \Scr K$. An easy calculation shows that 
a unitary operator  $u$ is a $K$-skeleton unitary if and 
only 
if   $a$ is a Fredholm partial isometry  on the submodule 
$p\Scr H_{\Scr A}$ 
and $d$ is a Fredholm partial isometry on the submodule 
$(1-p)\Scr H_{\Scr A}$; in other words,
all  $p-aa^*, p-a^*a, (1-p)-dd^*,
(1-p)-d^*d$  are   projections in $\Scr A\otimes \Scr K$. 
The term 
`{\it K}-skeleton' is chosen, since $K_0(\Scr A)$ is 
completely 
described by the homotopy classes of all such unitaries.

Let $ {GL^p_r(\Scr A)}$ be the topological group 
consisting of all those 
invertible operators in $\Scr L(\Scr H_{\Scr A})$ such that 
$xp-px\in \Scr A\otimes \Scr K$, equipped with the norm 
topology from 
$\Scr L(\Scr H_{\Scr A})$. Let $ {GL^p_{\infty}(\Scr A)}$ 
stand for the 
path component of $GL^p_r(\Scr A)$ containing the identity; 
in the special case when $\Scr A=\bold{C}$, we instead use 
the notation
$ {GL^p_r(\Scr H)}$ and $ {GL^p_{\infty}(\Scr H)}$, 
respectively.
Let $ {GL_{\infty }(\Scr A)}$ and $ {GL^0_{\infty}(\Scr 
A)}$ 
denote the group of all invertible elements in the 
unitization of 
$\Scr A\otimes \Scr K$ and its identity path component, 
respectively.

\proclaim {\num{2.1.} $K$-skeleton factorization theorem  
\cite{Zh4}\rm}
{\rm (i)  }  If $x\in GL^p_r(\Scr A)$, then  there exist
an element $k\in \Scr A\otimes \Scr K$, an invertible 
element 
$(\smallmatrix z_1&0\\0&z_2\endsmallmatrix)$, and a 
K-skeleton unitary $(\smallmatrix 
a&b\\c&d\endsmallmatrix)$ along 
$p$ such that $1+k\in GL^0_{\infty}(\Scr A)$ and 
$$x=(1+k)\left (\matrix z_1&0\\0&z_2\endmatrix \right )
\left (\matrix a&b\\c&d\endmatrix \right ).$$
A factorization of $x$ with the form above is called
\it a $K$-skeleton factorization  along $p$.\rm 

{\rm (ii) } If two $K$-skeleton factorizations of $x$ 
along $p$ are given,
say 
$$x=x_0x_p\left (\matrix a&b\\c&d\endmatrix \right 
)={x}'_o{x}'_p
\left (\matrix {a}'&{b}'\\{c}'&{d}'\endmatrix \right ),$$ 
then $[cc^*]-[bb^*]=[{c}'c^{\prime *}]-[{b}'b^{\prime *}] 
\in  K_0(\Scr A);$
in other words, $[cc^*]-[bb^*]$ is an invariant  
independent of  all 
(infinitely many) possible K-skeleton factorizations of 
$x$ along $p$.
\endproclaim

\demo{Outline of a proof} There is a shorter proof solely 
for this theorem.
For the sake of clarifying some internal relations among 
$\pi _0(GL^p_r(\Scr A))$, $ \pi _0([p]_0)$, and $K_0(\Scr 
A)$, 
we outline  a proof as follows.
First,  every  element in $GL^p_{\infty}(\Scr A)$ can be 
written as a product of the form 
$x_0x_p$ for some invertible $x_0\in GL^0_{\infty}(\Scr A)$
 with $x_0-1\in \Scr A\otimes \Scr K$ 
and another invertible $x_p$ with $x_pp=px_p$ [Zh4].
Secondly, write the  polar decomposition 
$x=(xx^*)^{{1}/2}u $, where 
$(xx^*)^{{1}/2}\in GL^p_{\infty}(\Scr A)$ and $u$ is a 
unitary in 
$GL^p_r(\Scr A)$.
Then consider the following subsets of
  $Gr^{\infty }(\Scr A)$:
$$\align 
{[upu^*]_r}&:=\{wupu^*w^*: 
w\in GL^0_{\infty }(\Scr A)\ \ \text{with}\ \ 
ww^*=w^*w=1\}\\
\intertext{and}
{[p]_0}&:=\{vpv^*: v\in GL^p_r(\Scr A)\ \ \ 
vv^*=v^*v=1\}\. \endalign$$
 Technical arguments show that $[upu^*]_r$ is precisely 
the path component of $[p]_0$ containing $upu^*$. 
Thirdly, there is a representative in $[upu^*]_r$ with the 
form 
$(p-r_1)\oplus r_2$ for some projections  
$r_1, r_2 \in \Scr A\otimes \Scr K$.  It follows that 
there exists a unitary $u_0\in GL^0_{\infty }(\Scr A)$ 
such that $$u_0^*upu^*u_0=(p-r_1)\oplus r_2.$$
Then one obtains a $K$-skeleton unitary 
$(\smallmatrix a&b\\c&d\endsmallmatrix )$ such that 
$u=u_0 (\smallmatrix a&b\\c&d\endsmallmatrix ),$
 where $bb^*=r_1$ and  $cc^*=r_2$. Since 
$(xx^*)^{{1}/2}u_0\in 
GL^p_{\infty}(\Scr A)$, we can rewrite it as a product in 
the desired form 
$x_0 (\smallmatrix z_1&0\\0&z_2\endsmallmatrix )$. 
The details are contained  in [Zh4].
\enddemo

  It follows from Theorem 2.1 that 
$x\cdot GL^p_{\infty}(\Scr A)= (\smallmatrix 
a&b\\c&d\endsmallmatrix)
\cdot GL^p_{\infty}(\Scr A)$ (cosets)  for each $x\in 
GL^p_r(\Scr A)$.
The invariant $[cc^*]-[bb^*]$ associated with the 
$K$-skeleton 
factorization of $x\in GL^p_r(\Scr A)$ 
yields the bijection $$\left (\matrix a&b\\c&d\endmatrix 
\right )\cdot 
GL^p_{\infty}(\Scr A)\ \longleftrightarrow \ 
[(p-bb^*)\oplus (cc^*)]_r.$$
It can be shown that 
$[(p-r_1)\oplus {r}'_1]_r=[(p-r_2)\oplus {r}'_2]_r$ 
iff $[{r}'_1]-[r_1] =[{r}'_2]-[r_2]$ in $K_0(\Scr A).$
Therefore, we  conclude the following theorem whose 
details are given in 
[Zh4].

\proclaim {\num{2.2.} Theorem  \cite{Zh4}\rm}
The maps  defined by 
$$
\pmatrix a & b\\ c & d\endpmatrix\cdot GL^p_{\infty}(\Scr 
A)\ \longmapsto\ 
[(p-r_1)\oplus r_2]_r\ \longmapsto \ [r_2]-[r_1]$$
are two bijections, which 
induce the following isomorphisms\RM:
$$GL^p_r(\Scr A)/{GL^p_{\infty }(\Scr A)}\ \cong\ D_h([p]_0)
\ \cong \ K_0(\Scr A),$$ 
where $GL^p_r(\Scr A)/{GL^p_{\infty }(\Scr A)}$ is the 
quotient 
group  with the induced  multiplication and 
$$D_h([p]_0)=\{[upu^*]_r: u\in GL^p_r(\Scr A)\ \ 
\text{with}\ 
uu^*=u^*u=1\}$$ is the set of all path components of 
$[p]_0$. 
The group operation on $D_h([p]_0)$ is defined by
$$[(p-r_1)\oplus {r}'_1]_r+[(p-r_2)\oplus {r}'_2]_r
=[(p-r_1-s_2)\oplus ({r}'_1\oplus {s}'_2)]_r$$
for some projections $s_2\in p(\Scr A\otimes \Scr K)p$  and 
${s}'_2\in  (1-p)(\Scr A\otimes \Scr K)(1-p)$ such that
$s_2\sim r_2, $ $s_2r_1=0,$  ${s}'_2\sim {r}'_2,$ and 
${s}'_2{r}'_1=0$.
\endproclaim

\proclaim {\num{2.3.} Theorem}  Let the base point of 
$[p]_0$ be $p$
and the base point of $GL^p_r(\Scr A)$ be the identity. 
Then\endproclaim
$$\align
 &\pi _{2k+1}([p]_0) \cong  \pi _{2k+1}(GL^p_r(\Scr A)) 
\cong 
K_1(\Scr A),\\
\intertext{and}
&\pi _{2k+2}([p]_0) \cong  \pi _{2k+2}(GL^p_r(\Scr A)) 
\cong K_0(\Scr A)
\ \ \ \forall\  k\ge 0.\endalign$$

\demo{Outline of a proof } Let $ {U_{\infty }(\Scr A)}$ be 
the 
unitary group of the unitization of $\Scr A\otimes \Scr K$,
and let $ {U_p(\Scr A)}$ be the subgroup of $U_{\infty 
}(\Scr A)$
consisting of all those unitaries commuting with $p$. 
First, the map $ \psi _p :  U_{\infty }(\Scr A) 
\longrightarrow  [p]_r$
defined by $\psi _p(u)= upu^*$ is a Serre (weak) fibration 
with 
a standard fiber $U_p(\Scr A)$ [Zh6, \S 2].  Secondly,  the 
long exact sequence of homotopy groups associated with 
this fibration 
breaks into short exact sequences [Zh6, 2.5, 2.8]:
$$0\longrightarrow \pi _{k+1}([p]_r)\longrightarrow 
\pi _k(U_p(\Scr A) )\longrightarrow 
\pi _k(U_{\infty }(\Scr A))\longrightarrow 0\qquad (k\ge 
0).$$
Thirdly, by an analysis on this  short exact sequence one 
concludes 
$$\pi _{2k+2}([p]_0)\cong K_0(\Scr A)\ \ \ \text{and}\ \ \ 
\pi _{2k+1}([p]_0)\cong K_1(\Scr A)\qquad(k\ge 0).$$

It is well known that the subgroup $ {U^p_r(\Scr A)}$ 
consisting of all unitary elements in $GL^p_r(\Scr A)$ is 
homotopy equivalent to $GL^p_r(\Scr A)$. We consider the 
maps 
$U^p_r(\Scr A) \longrightarrow  [p]_0$ defined by $\phi 
_p(u)=upu^*$.
It can be shown that
$\phi _p$ is a weak fibration with a standard fiber $ 
{U^p(\Scr A)}$,
where $ {U^p(\Scr A)}$ is the group consisting of all 
those unitaries in $U^p_r(\Scr A)$ commuting with $p$.
An argument similar to that above applies to this fibration.
One can show that 
$\pi _{2k+1}(U^p_r(\Scr A))\cong K_1(\Scr A)$ and
$\pi _{2k+2}(U^p_r(\Scr A))\cong K_0(\Scr A)$  for $k\ge 0.$
The details are given in [Zh6, \S 4]. \enddemo

\subheading{\num{2.4.} Special case ${\Scr A=C(X)}$\rm}
In particular, if $\Scr A$ is taken to be  the commutative 
C*-algebra 
$C(X)$ consisting of all complex-valued continuous 
functions on a compact 
Hausdorff space $X$, then each element in $\Scr L(\Scr 
H_{C(X)})$ 
can be identified with a norm-bounded, *-strong continuous 
map 
from $X$ to $\Scr L(\Scr H)$ [APT]. 
Here $\Scr L(\Scr H)\supset \{x_{\lambda }\}$ converges to 
$x$ 
in the *-strong operator topology iff 
$$\Vert (x_{\lambda }-x)k\Vert +\Vert k(x_{\lambda 
}-x)\Vert \rightarrow 0
\quad \text{for any } k\in \Scr K.$$
Obviously, $\Scr L(\Scr H_{C(X)})$ contains the C*-tensor 
product 
$\Scr L(\Scr H)\otimes C(X)$ 
consisting of all norm-continuous maps from $X$ to $\Scr 
L(\Scr H)$
as a C*-subalgebra. Then Theorems 2.1 and 2.2 in this 
special case
are interpreted as follows.

\proclaim {\num{2.5.}  Corollary } Let $ {GL^{\infty 
}(\Scr H)}$ be the 
group of all invertible operators in $\Scr L(\Scr H)$.

{\rm (i)   } If $f:X\longrightarrow GL^{\infty }(\Scr H)$ 
is a 
norm-bounded, *-strong continuous map and $p$ is a 
projection in 
the infinite Grassmann space $Gr^{\infty }(\scr H)$ such 
that 
$pf-fp\in \Scr K\otimes C(X)$, then $f$ can be factored as 
the following  
product of three invertible maps
$$f(.)=\left (\matrix 1+k_{11}(.)&k_{12}(.)\\k_{21}(.)&1+
k_{22}(.)
\endmatrix \right )\left (\matrix 
g_1(.)&0\\0&g_2(.)\endmatrix \right )
\left (\matrix a(.)&b(.)\\c(.)&d(.)\endmatrix \right );$$
where $k_{ij}(.)$\RM's are norm-continuous maps from $X$ 
to $\Scr K$, 
$g_1(.)\oplus g_2(.)$ is a  norm-bounded, *-strong 
continuous map
from $X$ to $GL^{\infty }(\Scr H)$, $a(.)$, $d(.)$ are 
*-strong continuous maps from $X$ to the set of Fredholm 
partial 
isometries on $p\Scr H$ and $(1-p)\Scr H$, respectively, 
and $c(.)$, 
$b(.)$ are norm-continuous maps from $X$ to the set of 
partial isometries 
in $\Scr K$. Furthermore, $$[c(.)c(.)^*]-[b(.)b(.)^*]\ \in 
K_0(C(X))\ (\cong K^0(X))$$ is an invariant independent of 
all 
possible factorization with  the above form. 

{\rm (ii)  } The groups $[X, GL^p_r(\Scr H)]$, $[X, [p]_0]$,
and $K_0(C(X))$ are isomorphic,  where $[X,.]$ is the set 
of 
homotopy classes of norm-bounded, *-strong continuous maps 
from $X$ to $(.)$. 
\endproclaim

\subheading{\num{2.6.} Invertible dilations of a Fredholm 
operator  }
Let us illustrate  a $K$-skeleton factorization of  any 
invertible dilation of a Fredholm operator $x\in \Scr 
L(\Scr H_{\Scr A})$. 
There are of course infinitely many  invertible $2\times 
2$ matrices 
with the form $$D_2(x):=\left (\matrix 
x&y_1\\y_2&z\endmatrix \right )\ \in 
M_2(\Scr L(\Scr H_{\Scr A})).$$
Each such $2\times 2$ invertible matrix  is called 
\it an invertible dilation of $x$. \rm 
Specific constructions of such a dilation were 
given by P. Halmos [{Ho}, 222] and A.  Connes [{Co}].
For each invertible dilation of $x$ it follows from the  
$K$-skeleton Factorization Theorem 2.1 that
$$\left (\matrix x&y_1\\y_2&z\endmatrix \right )=
\left (\matrix 1+a_{11}&a_{12}\\a_{21}&1+a_{22}\endmatrix 
\right )
\left (\matrix z_1&0\\0&z_2\endmatrix \right )
\left (\matrix v&1-vv^*\\1-v^*v&-v^*\endmatrix \right ),$$
where $a_{ij}$'s are  some elements in  $\Scr A\otimes 
\Scr K$, $z_1, z_2
\in GL^{\infty }(\Scr A)$, and the above matrix  on the 
right, say $w$, 
is a familiar unitary matrix occurring in the index map in 
$K$-theory [{Bl}, 8.3.2] in which $v$ is a Fredholm 
partial isometry in 
$\Scr L(\Scr H_{\Scr A})$.  Set $p=\text{diag}(1, 0)$.
It is well known that $$[1-v^*v]-[1-vv^*]\ \in K_0(\Scr A)$$
is precisely the Fredholm index 
$\text{Ind}(v)=\text{Ind}(pxp)$\ (on
$p\scr H_A)$.
It follows from Theorem 2.1(ii) that those $K$-skeleton 
unitaries
associated with all possible invertible dilations of $x$ 
in $M_2(\Scr L(\Scr H_{\Scr A}))$ only differ from $w$ by 
a factor 
in $GL^p_{\infty}(\Scr A)$. 

\heading 3. Factorizations of invertible operators 
with integer indices \endheading

Now we consider some special  cases such that $K_0(\Scr 
A)\cong Z$ 
(the group of all integers); for example, $\Scr 
A=\bold{C}$, 
or $\Scr A=C(S^{2n+1})$ where $S^m$ is the standard 
$m$-sphere, 
or $\Scr A=\Scr O_{\infty}$, the Cuntz algebra generated 
by isometries 
$\{s_i\}_{i=1}^\infty \subset \Scr L(\Scr H)$ such that 
$\sum_{i=1}^\infty s_is_i^*\le 1$.

Let $p$ be any projection in $Gr^{\infty }(\Scr H)\subset 
Gr^{\infty }(\Scr A)$ [the inclusion holds because $\Scr 
H\subset 
\Scr H_{\Scr A}$ and $\Scr L(\Scr H)\subset \Scr L(\Scr 
H_{\Scr A})$].
Let $\{\xi _i\}_{i=0}^{+\infty }$ be any orthonormal basis 
of
the subspace $p\Scr H$ and $\{\xi_i\}^{i=-1}_{-\infty}$ be 
any orthonormal 
basis of the subspace $(1-p)\Scr H$. Then $\{\xi 
_i\}_{-\infty }^{+\infty }$
is an orthonormal basis of both $\Scr H$ and $\Scr H_{\Scr 
A}$.
Let $u_0$ denote the bilateral shift associated with 
the basis $\{\xi _i\}_{-\infty }^{+\infty}$ of $\Scr H$,
defined by  $u_0(\xi _i)=\xi _{i+1} $ for all $i\in Z.$
Clearly, $u_0$ is a $K$-skeleton unitary of $\Scr L(\Scr 
H_{\Scr A})$ along $p$.
Applying the $K$-skeleton Factorization Theorem 2.1 to the 
above special cases, 
we have the following factorizations of invertible operators
orientated by the integer-valued Fredholm index:

\proclaim {\num{3.1.}  Corollary }   Suppose that 
$K_0(\Scr A)\cong Z$ is generated by $[1]$ where $1$ is the 
identity of $\Scr A$. If $x$ is  an invertible operator on 
$\Scr H_{\Scr A}$ 
such that $px-xp\in \Scr A\otimes \Scr K$, then   $x=(1+
k)x_pu_0^{-n},$
where   $k\in \Scr A\otimes \Scr K$,  $x_p$ is an 
invertible 
operator commuting with $p$, and the integer $n$  is the 
Fredholm index 
  of $pxp$ on the submodule $p\Scr H_{\Scr A}$, say
$\roman{Ind}(pxp)$, which is  independent of the choice of 
$\{\xi _i\}_{i=0}^{+\infty }$, $\{\xi _i\}^{-1}_{-\infty 
}$ and 
all possible factorizations along $p$ with the same form 
above. \endproclaim 
\demo{Outline of a proof} 
It is obvious that $\text{Ind}(pu^n_0p)=-n$.
Let $G$ be the group $\{u_0^n: n\in Z\}$ in which every 
element is 
a $K$-skeleton unitary along $p$. As a special case of 
Theorem 2.1 one can show that the map from $G$ to
$GL^p_r(\Scr A)/
{GL^p_{\infty }(\Scr A)}$ defined by 
$u_0^n\longmapsto u_0^n\cdot GL^p_{\infty}(\Scr A)$
is a group  isomorphism. It follows that 
$\pi _0(GL^p_r(\Scr A))=\{u_0^n\cdot GL^p_{\infty }(\Scr 
A): n\in Z\}.$
Then the factorization follows.  The reader may want to 
consider 
the extreme case $\Scr A=\bold{C}$ and then generalize 
the conclusion to a larger  class of C*-algebras.
\enddemo 

A similar proof yields the following alternative  
factorization
of $x$ as a product of three invertibles 
under the same assumptions as of Corollary 3.1:
$$x=\cases (1+k_1)x_1 & \  \text{if}\ \text{Ind}(pxp)=0,\\
(1+k_2)x_2(u_1\oplus u_2\oplus\dotsb\oplus u_{-n}\oplus 
w_1) & \ \text{if}
\ \text{Ind}(pxp)=n<0,\\
(1+k_3)x_3(u^*_1\oplus u^*_2\oplus\dotsb\oplus u^*_n\oplus 
w_2) & \ \text{if}
\ \text{Ind}(pxp)=n>0,\endcases $$
where $u_i$ is a bilateral shift on a subspace $\Scr H_i$ of
$\Scr H$ for $1\le i\le n$, $w_j$\<'s are unitary 
operators on 
$(\bigoplus _{i=1}^n\Scr H_i)^{\bot}$,  $k_j\in \Scr 
A\otimes \Scr K$,
and $x_j$'s are invertible operators commuting with $p$.

\proclaim {\num{3.2.} Corollary  } Suppose that
$K_0(\Scr A)\cong Z$ is generated by $[1]$. 
If $x$ is an \it arbitrary \rm     element $\Scr L(\Scr 
H_{\Scr A})$ and 
$p\in Gr^{\infty }(\Scr H)$ (as above) such that $px-xp\in 
\Scr A\otimes 
\Scr K$, then there exists a \it unique \rm 
norm-continuous map 
$x(\lambda )$ from  $C\setminus \sigma (x)$ to 
$GL^p_{\infty}(\Scr A)$, 
where $\sigma (x)$ is the spectrum of $x$,  such that 
$x-\lambda =x(\lambda )u_0^{-n_i},$ where  
$n_i=\text{Ind}(p(x-\lambda _i)p)$
and $\lambda _i$ is any complex number in the $i$\<th path 
component 
$O_i$ of $C\setminus \sigma (x)$. 
An alternative K-skeleton factorization of $x-\lambda $ 
for $\lambda \in O_i$ is as follows (when $n_i\neq 0)$:
$$x-\lambda =\cases y_i(\lambda )(u_1\oplus u_2\oplus 
\dotsb\oplus u_{\vert 
n_i\vert }\oplus w_i)&\  \text{if}\ \ 
\text{Ind}(p(x-\lambda _i)p)=n_i<0,\\
{y}'_i(\lambda )(u^*_1\oplus u^*_2\oplus \dotsb\oplus 
u^*_{n_i }\oplus v_i) 
&\  \text{if}\ \ \text{Ind}(p(x-\lambda _i)p)=n_i> 0,
\endcases$$
where $u_i$\<'s are bilateral shifts on mutually 
orthogonal closed subspaces 
$\Scr H_i$'s of $\Scr H$,  $w_i, v_i$\<'s  are  unitary 
operators on 
the subspace $(\bigoplus _{i=1}^{\vert n_i\vert} \Scr 
H_i)^{\bot }$\<, and
$y_i(\lambda ),{y}'_i(\lambda )$
are norm-continuous maps from $O_i$ to
$GL^p_{\infty}(\Scr A)$.
\endproclaim 

\subheading{\num{3.3.}  Winding numbers of invertible 
operators }
Using the first factorization in Corollary 3.2, 
we assign an integer $n_i$ to each path component $O_i$ of
$C\setminus \sigma (x)$,
which is precisely the minus winding number of $u_0^{-n_i}$
as a continuous map from $S^1$ to $S^1$
(via the Gel'fand transformation).  We call $n_i$  the \it 
winding number
of $x$ along  $p$ over $O_i$.   \rm 
As a particular case,  if $x$ is an operator whose 
essential 
spectrum, the spectrum of $\pi (x)$ in the generalized 
Calkin 
algebra $ \Scr L(\Scr H_{\Scr A})/{\Scr K(\Scr H_{\Scr 
A})}$,
 does not separate the plane, then all  winding numbers of 
 $x$ along 
any $p\in Gr^{\infty }(\Scr A)$  are  zero 
as long as  $px-xp\in \Scr A\otimes \Scr K$. 
There is another way to describe the integer $n_i$.

\proclaim {\num{3.4.}  Corollary  }
Let  $G_i(x)$ denote the subgroup  of $GL^p_r(\Scr A)$ 
generated by
 $GL^p_{\infty}(\Scr A)$ and $x-\lambda _i$  where 
$\lambda _i\in O_i$. Then $G_i(x)/{GL^p_{\infty}(\Scr 
A)}\cong n_iZ$,
and hence 
$GL^p_r(\Scr A)/{G_i(x)} \cong Z_{n_i},$ the finite
cyclic group of order $n_i.$\endproclaim 
 In particular, one can apply the above factorizations 
to an invertible dilation of a pseudodifferential operator 
of order zero 
on a compact manifold and classical multiplication 
operators. 
Let us spend  few lines  to look at the following familiar 
examples.

\subheading{\num{3.5.} Multiplication operators  }
Let $M_f$ be the invertible multiplication operator with  
symbol $f$ in $L^{\infty}(S^1)$, where $S^1$ is the unit 
circle; i.e., 
$M_f(g) =fg$ for any $ g\in L^2(S^1).$ If  
$p$ is a projection on $L^2(S^1)$  such that 
$\text{dim}(p)=\text{codim}(1-p)=\infty $ and 
$pM_f-M_fp$ is a compact operator,
then it follows from Corollary 3.1 that 
$M_f=(1+k)x_pu_0^{-n},$ where 
$n=\text{Ind}(pM_fp)$, $k$ is a compact operator on 
$L^2(S^1)$,
 $x_p$ is an invertible operator on $L^2(S^1)$ commuting 
with $p$,
and $u_0$ is a bilateral shift operator associated with a
fixed orthonormal basis of $L^2(S^1)$. 
It is well known that $pM_fp$ is  a familiar Toeplitz 
operator on 
the subspace  $pL^2(S^1)$.

\subheading{\num{3.6.} Restricted loop group along ${p}\in 
Gr^{\infty }(\Scr H)$ }
Consider the following \it restricted loop group along $p$
\rm  consisting of all norm-bounded, 
*-strong continuous maps from $S^1$ to $GL^p_r(\Scr H)$,
denoted by $\roman{Map}(S^1,GL^p_r(\Scr H))_{\beta}.$ 
Since
 $K_0(C(S^1))=Z$, each $f\in \text{Map}(S^1,GL^p_r(\Scr 
H))_{\beta }$
can be factored as 
$f=
(1+f_0)f_1u_0^{-n},$
where $n=\text{Ind}(pfp)$, $f_0$ is a norm-continuous map 
from 
$S^1$ to $\Scr K$,  $f_1$ is a
 *-strong continuous map from $S^1$ 
to $GL^{\infty}(\Scr H)$ such that
$f_1(z)p=pf_1(z)$ for any $z\in S^1$, and $u_0$ is a 
bilateral shift 
with respect to a fixed orthonormal basis of $\Scr H$.
If $f$ is norm-continuous, then $f_1$ is also norm 
continuous.
Furthermore, $[S^1, GL^p_r(\Scr H)] \cong [X, [p]_0] \cong 
 Z.$
The same conclusions also  hold, if $S^1$ is replaced by 
$S^{2n+1}$ for any
$n\geq 1$.

\rem{\num{3.7.} Remarks}
(i)
Theorems 2.1--2.3 still hold, if $\Scr A$ 
is any stably unital C*-algebra; i.e., $\Scr A\otimes \Scr 
K$ has 
an approximate identity consisting of a sequence of 
projections 
[Bl, 5.5.4; Zh4]. 

(ii) Let $\text{Index}(x,p)$ denote the invariant 
$[cc^*]-[bb^*]\in K_0(\Scr A)$ in Theorem 2.1(ii). If $p$ 
is fixed, 
then $\text{Index}(x,p)$ is precisely the Fredholm index  
of $pxp$ 
as an operator on $p\Scr H_{\Scr A}$ and fits into 
the established theory of the $K_0(\Scr A)$-valued 
Fredholm index. 
However, some new results do arise from invariants of 
$\text{Index}(x,p)$ as the variable  $p$  runs in 
$\{p\in Gr^{\infty }(\Scr A): xp-px\in \Scr A\otimes \Scr 
K\}$ 
or as $x$ and $p$ jointly change [Zh7]. 
As a matter of fact, $\text{Index}(x,p)$ is an invariant 
under homotopy 
and perturbation by elements in $\Scr A\otimes \Scr K$ 
with respect to 
both variables $x$ and $p$.
For example, by the combination of the $K$-skeleton  
Factorization Theorem and 
certain invariants of $\text{Index}(x,p)$,  we  proved 
[Zh7] the following: 

\proclaim{Theorem}
$$\pi _0(GL({M_n(\bold {C})}'_e)) \ \cong \ 
\{k\in K_0(\Scr A):n\cdot k=0\}\quad\text{ for any $n\geq 
2$};$$
where $GL({M_n(\bold{C})}'_e)$ denotes  the group of all 
invertibles 
in the essential commutant ${M_n(\bold{C})}'_e$ of 
$M_n(\bold{C})$ which is naturally embedded in $M_n(\Scr 
L(\Scr H_{\Scr A}))$. 
\endproclaim

(iii) The reader may want to compare 
(3.1)--(3.3) and the famous BDF theory [BDF1,2] to see 
their obvious 
relations; we work with invertibles on $\Scr H_{\Scr A}$, 
while the BDF 
theory dealt with Fredholm operators.

(iv) In [PS]  Pressley and Segal have studied  
the \it restricted general linear group  
$${GL_{\roman{res}}(\Scr H)}
:=\{x\in GL^{\infty} (\Scr H): xp-px \ \text{is 
Hilbert-Schmidt }\}$$
and given some applications to  the Kdv equations. It is a 
hope 
that our results  will shed some light in the same 
direction.
\endrem

\enddocument